\documentclass[draft,a4paper,11pt,fleqn,leqno]{article}
\usepackage{amssymb,amsmath,cite}
\usepackage[latin1]{inputenc}
\usepackage{latexsym}

\numberwithin{equation}{section}
\newtheorem{theorem}{Theorem}[section]
\newtheorem{definition}{Definition}[section]
\newtheorem{lemma}{Lemma}[section]
\newtheorem{proposition}{Proposition}[section]
\newtheorem{remark}{Remark}[section]

\newtheorem{corollary}{Corollary}[section]
\def\proof{\text{\it Proof:} }

\newcommand{\ann}[1]{}

\def\fin{\mbox{$\blacksquare$}} 

\def\R{\mathbb{R}}

\begin{document}
\title{\bf Extinction of solutions of semilinear higher order parabolic equations with degenerate absorption potential}
\author{{\bf\large Y. Belaud, A. Shishkov}}\vspace{3mm}
\date{}
\maketitle

\begin{center}
{\bf\small Abstract}

\vspace{3mm} \hspace{.05in}\parbox{4.5in} {{\small We study the
first vanishing time for solutions of the Cauchy-Dirichlet\\
problem to the semilinear $2m$-order ($m \geq 1$) parabolic equation\\
$u_t+Lu+a(x) |u|^{q-1}u=0$, $0<q<1$ with $a(x) \geq 0$ bounded in
the bounded domain $\Omega \subset \R^N$. We prove that if $N>2m$
and $\displaystyle \int_0^1 s^{-1} \text{meas} \{x \in \Omega :
|a(x)| \leq s \}^\frac{2m}{N} ds < + \infty$, then the solution
$u$ vanishes in a finite time. When $N=2m$, the condition becomes
$\displaystyle \int_0^1 \!\!\!\!s^{-1} \left( \text{meas} \{x \in
\Omega : |a(x)| \leq s \}\right) \left( -\ln \text{meas} \{x \in
\Omega : |a(x)| \leq s \}\right) ds < + \infty$.}}
\end{center}
\vspace{0.5 cm} {\it \footnotesize Key words}. {\scriptsize
nonlinear equation, energy method, vanishing solutions, semi-classical analysis\\
MSC 35B40, 35K20, 35P15}
\section{Introduction and main results}

Let $\Omega\subseteq\mathbb R^N,\ N\geq1$, be arbitrary bounded
domain. In cylindrical domain $\Omega\times(0,\infty)$ we consider
 the following Cauchy-Dirichlet problem
\begin{eqnarray} \label{parabolicequation}
u_t + L(u) + a(x) f(u) = 0 \; \text{in} \; \Omega \times
(0,\infty), \ f(u)=f_1(u):=|u|^{q-1}u, \  0<q<1,
\end{eqnarray}
\begin{eqnarray}
D_x^\alpha u(x,t)=0 \; \text{on} \; \partial \Omega \times (0,\infty), \; \forall \alpha : |\alpha| \leq m-1,
\end{eqnarray}
\begin{eqnarray} \label{initialdata}
u(x,0)=u_0(x) \in L^2(\Omega).
\end{eqnarray}
Here $L$ is a divergent differential $2m$-order operator of the
form :
\begin{eqnarray} \label{operatorl}
L(u) = (-1)^m \sum_{|\alpha|=m} D_x^\alpha
a_{\alpha}(x,u,D_xu,\dots ,D^m_x u), \quad m\geqslant 1,
\end{eqnarray}
with Caratheodory functions $a_\alpha(x,\xi)$ (continuous with
respect $\xi$ and measurable with respect to $x$) satisfying
sublinear growth condition:
\begin{multline} \label{growth}
|a_{\alpha}(x,\xi)|\leqslant c \sum_{|\gamma|=m}|\xi_\gamma|\quad
\ \forall\,\xi=\{\xi_\gamma\}\in\mathbb R^{M(m)},\quad
|\alpha|\leqslant m, x\in\bar{\Omega};\,\, c=const,
\end{multline}
where $M(m)$ is the number of different multi-indices
$\gamma=(\gamma_1,\gamma_2,\dots \gamma_N)$ of the length
$|\gamma|:=\gamma_1+\dots +\gamma_N\leqslant m$, and the
absorptional potential $a(x)$ is nonnegative, measurable, bounded
function in $\Omega$.

Our main condition on the operator $L$ is the following coercivity
condition:
\begin{multline}\label{coercivity}
(L(v),v):=\int_\Omega
\sum_{|\alpha|=m}a_\alpha(x,v,\dots,D_{x}^{m}v)D^{\alpha}_{x}v\,dx\geqslant
C \int_\Omega \vert D^{m}_{x}v\vert^2\,dx \\ \forall v \in
W^{m,2}_{0}(\Omega), C=const>0,
\end{multline}
where $W^{m,2}_{0}(\Omega)$ is the closure in the norm
$W^{m,2}_{0}(\Omega)$ of the space $C_{0}^{m}(\Omega)$.
\begin{remark}\label{Rem1}
 Well known sufficient condition of
\eqref{coercivity} is
\begin{eqnarray}\label{monotonicity}
\sum_{|\alpha|=m}a_\alpha(x,\xi)\xi_\alpha\geqslant
C\sum_{|\gamma|=m}|\xi_\gamma|^2 \quad
\forall\,\xi=\{\xi_\beta\}\in\mathbb R^{M(m)}, \quad \forall\,x\in
\bar{\Omega}.
\end{eqnarray}
\end{remark}
\begin{remark}\label{Rem2}
 In the linear case for operator $L=\sum_{|\alpha|=|\beta|=m}D^\alpha
a_{\alpha\beta}D^\beta$ with constant coefficients
$a_{\alpha\beta}$ as it is easy to check by using of Plancherel
theorem, property \eqref{coercivity} is guaranteed by ellipticity
condition:
\begin{eqnarray}\label{ellipticity}
\sum_{|\alpha|=|\beta|=m}a_{\alpha\beta}\zeta^\alpha\zeta^\beta\geqslant
C|\zeta|^{2m} \quad
\forall\,\zeta=(\zeta_1,\dots\zeta_N)\in\mathbb
R^N,\zeta^\alpha:=\zeta^{\alpha_1}\zeta^{\alpha_2}\dots\zeta^{\alpha_N}.
\end{eqnarray}

\end{remark}

\begin{definition}
We will say that problem \eqref{parabolicequation}-\eqref{initialdata} has the extinction in finite time (EFT) property if for arbitrary solution $u$, there exists some positive $T_0$ such that $u(x,t)=0$ a.e. in $\Omega$, $\forall t \geq T_0$.
\end{definition}
Firstly EFT-property for simplest semilinear heat equation with
strong absorption was observed by A. S. Kalashnikov
\cite{Kalash1}. Later mentioned property (conditions of occurence
of extinction, estimates of extinction time, asymptotic of
solution near to the extinction time and so on) was investigated
for different classes of second order semilinear and quasilinear
parabolic equations of diffusion-absorption type by many authors
(see \cite{KnSt1,Pay,Ker2,Knerr1,BSt,FH,ChMM}). F. Bernis
\cite{Bern3} proved the EFT-property for energy solutions to
higher order semilinear  and quasilinear parabolic equations with
strong absorption. Dependence of extinction properties of energy
solutions to mentioned higher order equations on local structure
of initial function was studied in \cite{Shish1}. Extinction
properties for second order semilinear parabolic equations of
diffusion-absorption type with nondegenerate $(x,t)$-dependent
absorptional potential was studied in \cite{Kalash2,KerN,JL1,JL2}.

V. Kondratiev, L. Veron \cite{KV1} firstly initiated the study of
EFT-property for second order equation \eqref{parabolicequation}
($m=1$) in the case of degenerate absorptional potential $a(x)$:
\begin{equation}\label{1.1sh}
\qquad \inf\{a(x):x\in\Omega\}=0.
\end{equation}
It happens that occurence of mentioned property depends
essentially on the structure of the set of degeneration and on the
behaviour of potential $a(x)$ in the neighbourhood of this set.
They in \cite{KV1} considered homogeneous Neumann problem for
second order equation \eqref{parabolicequation} ($m=1$) and proved
the following general sufficient condition for EFT-property:
\begin{equation}\label{1.2sh}
\hskip120pt\qquad\sum_{i=1}^\infty\mu_i^{-1}<\infty,
\end{equation}
$$
\mu_k:=\inf\bigg\{\int_\Omega\sum_{i,j=1}^N(a_{ij}v_{x_i}v_{x_j}+2^ka(x)v^2)dx:v\in
W^{1,2}(\Omega),\ \int_\Omega v^2dx=1\bigg\},\quad
\forall\,k\geqslant 1.
$$
Method from \cite{KV1} (semiclassical or KV-method) was developed
in \cite{BHV01} and the following explicit sufficient condition of
EFT-property for Dirichlet and Neumann boundary problem for second
order equation \eqref{parabolicequation}  was established:
\begin{equation}\label{1.3sh}
\ln a(x)^{-1}\in L^p(\Omega)\text{ for some }p>\frac N2.
\end{equation}
As a consequence, if $\{0\}\in\Omega$, then arbitrary potential
\begin{equation}\label{1.4sh}
a(x):a(x)\geq a_\alpha(|x|):=\exp(-\frac1{|x|^\alpha})\ \forall
x\in\Omega
\end{equation}
satisfies condition \eqref{1.3sh} by arbitrary $\alpha<2$. On the
other hand, for potential $a_\alpha(x)$ with $\alpha>2$
EFT-property fails \cite{BHV01}.

In \cite{BS07} there was elaborated the adaptation of local energy
method from \cite{KSh,Shish1}  to the study of extinction
properties of energy solution to second order parabolic equations
with radial degenerate absorptional potential. As result the
following sharp Dini-like sufficient condition of EFT-property was
obtained:
\begin{equation}\label{1.5sh}
a(x)\geq\exp(-\frac{\omega(|x|)}{|x|^2}),\text{ where }
\omega(s)>0\,\forall\,s>0:\omega(0)=0,\,\int_0^1\frac{\omega(s)}sds<\infty.
\end{equation}

The drawback of using regularizing effects does not enable the
KV-method to be extended to higher order operators but for small
dimensions (continuous injection of $W^{m,2}(\Omega)$ into
$L^\infty(\Omega)$). Moreover, the the local energy estimate
method from \cite{BS07} is developed till now for radial potentials
$a(|x|)$ only. These reasons lead us to construct some new variant
of semiclassical method. On the contrary with \cite{KV1}, we
consider a family of first eigenvalues of non-linear Schr\"odinger
operator directly connected with equation
\eqref{parabolicequation}, instead of eigenvalues $\mu_i$
\eqref{1.2sh} of auxiliary linear Schr\"odinger operator. As a
consequence, we do not need regularizing effects for solutions of
problem \eqref{parabolicequation}-\eqref{initialdata}. But this
means also that we can not use Lieb-Thirring formula \cite{LT76}
to estimate the first eigenvalue. Therefore, we provide
estimations of eigenvalues thanks to suitable Sobolev embedding
inequalities.

Thus let us denote for arbitrary potential $a(x)\geqslant 0$ the
function
\begin{equation}\label{1.6sh}
M_a(s):=\text{meas}\{x\in\Omega:a(x)\leq s\}.
\end{equation}
Then main assumptions on the degeneration of $a(x)$ are:
\begin{eqnarray} \label{integralassumption}
\int_0^1 s^{-1}M_a(s)^\theta ds < + \infty, \; \theta = \min\left(\frac{2m}{N},1\right), N \neq 2m,
\end{eqnarray}
and
\begin{eqnarray} \label{integralassumptionorlitz}
\int_0^1 s^{-1}M_a(s) \left(-\ln M_a(s)\right) ds < + \infty \; \text{for} \; N=2m.
\end{eqnarray}
For a function $a(x)$ satisfying \eqref{integralassumption} or
\eqref{integralassumptionorlitz} the set where it takes small
values is small enough. For instance, if $a(x)\geq\gamma>0$ then
$M_a(s)=0\quad\forall\,s<\gamma$ and, as consequence, integrals
are finite. On the contrary, if $a(x)=0$ on a set of positive
measure, integrals becomes infinite.
\begin{definition} \label{definitionregularity}
A function $u$ belonging to the space $\displaystyle {\cal
C}([0,\infty),L^2(\Omega)) \bigcap L^2_{loc}([0,\infty),\break
W^{m,2}_0(\Omega))$ is a weak solution of problem
\eqref{parabolicequation}-\eqref{initialdata} if initial condition
\eqref{initialdata} holds and if for any $\zeta \in
L^2_{loc}([0,+\infty),W^{m,2}_0(\Omega))$, there holds
\begin{eqnarray} \label{weaksolutionparabolic}
\int_0^T \!\left[ \!\langle u_t(t,.),\zeta(t,.) \rangle +
\!\int_\Omega \!\left( \sum_{|\alpha|=m} \!\!\!\!\!\!
\,\,\,\,\,a_{\alpha}(x,u,\dots ,D^m_x u)D_x^\alpha \zeta + a(x)
|u|^{q-1} u \zeta \right) \!dx \right] \!dt = 0,
\end{eqnarray}
for all $T>0$, where $\langle .,. \rangle$ is the paring of
elements from $(W^{m,2}_0(\Omega))^*$ and $W^{m,2}_0(\Omega)$,
$W^{m,2}_0(\Omega)$ being the closure in the norm
$W^{m,2}(\Omega)$ of the space ${\cal C}_0^m(\Omega)$.
\end{definition}
Main results are the following :
\begin{theorem} \label{theorem}
Under assumptions \eqref{operatorl} and \eqref{monotonicity},
\begin{itemize}
    \item[a)] if $N \neq 2m$ and \eqref{integralassumption} holds, then all weak solutions of problem \eqref{parabolicequation}-\eqref{initialdata} have the EFT-property,
    \item[b)] if $N = 2m$ and \eqref{integralassumptionorlitz} holds, then all solutions of problem \eqref{parabolicequation}-\eqref{initialdata} have the EFT-property.
\end{itemize}
\end{theorem}
\begin{remark}\label{r1.1}
Condition \eqref{integralassumption} by $m=1,N>2$ in the case of radial potential $a(|x|)$ implies condition \eqref{1.5sh}.
\end{remark}
Additionally, for second order equation ($m=1$), we improved the
KV-method for the Dirichlet problem and the results holds as
follows:
\begin{theorem}\label{theorem1.2}
Let $m=1$ in equation \eqref{parabolicequation} and
\begin{eqnarray} \label{integralassumptionlaplacian}
\int_0^1 s^{-1}M_a(s)^\frac{2}{N} ds < + \infty,
\end{eqnarray}
then all solutions of problem \eqref{parabolicequation}-\eqref{initialdata} have the EFT-property.
\end{theorem}

\section{Proof of Theorem \ref{theorem}}
The general principle is to find a lower bound for the function\\
$U(t):=\displaystyle \int_\Omega (|D_{x}^{m}u|^2 + a(x) |u|^{1+q})
dx$, with the help of the function
\begin{equation} \label{lambdaone}
\lambda_1(h) = \inf \left\{\int_\Omega (|D_x^m v|^2 + a(x)
|v|^{1+q})dx, v \in W^{m,2}_0(\Omega),||v||^2_{L^2(\Omega)}=h
\right\}
\end{equation}
 The key-stone of this section is the following :
\begin{proposition} \label{propositionkeystone}
If
\begin{equation}
\int_0^1 \frac{1}{\lambda_1(h)} dh < + \infty,
\end{equation}
then all  solutions of problem
\eqref{parabolicequation}-\eqref{initialdata} vanish in a finite
time and in this case,
\begin{equation}
T \leq \frac{1}{2} \int_0^{||u_0||^2_{L^2(\Omega)}} \frac{1}{\lambda_1(h)} dh .
\end{equation}
\end{proposition}
\proof Using $\zeta=u$ in \eqref{weaksolutionparabolic} gives for
all $0 \leq s < t$,
\[
\int_s^t \left[ \langle u_\tau(\tau,.),u(\tau,.) \rangle +
\!\int_\Omega \!\left( \sum_{|\alpha|=m} \!\!\!\!\!\!\!\,\,\,
a_{\alpha}(x,u,\dots,D^{m}_{x}u)D_x^\alpha u(x,\tau)  +a(x)
|u|^{q+1} \right) dx \right] \!d\tau = 0,
\]
which implies  by formula of integration by parts (see
\cite{Bern3}),
\[
\!\!\!\!\!\!\!\frac{1}{2}\!\int_\Omega \!\!\left( |u(t,.)|^2 -
|u(0,.)|^2\right) dx + \!\int_0^t \!\! \!\int_\Omega \!\left(
\sum_{|\alpha|=m} \!\!\!\!\!\!\! \quad a_\alpha(x,u,\dots ,
D^m_xu) D_x^\alpha u + a(x) |u|^{q+1} \right) dx \, \!d\tau = 0.
\]
But the second term is absolutely continuous with respect to time. Therefore the first term is also absolutely continous and has derivative a.e. with respect to time which leads to
\[
\frac{1}{2} \frac{d}{dt} (||u||^2_{L^2(\Omega)}) + \int_\Omega
\!\left( \sum_{|\alpha|=m} \!\!\!\!\!\!\!\quad a_\alpha(x,u,\dots
, D^m_xu) D_x^\alpha u + a(x) |u|^{q+1} \right) dx = 0.
\]
Clearly, from the property \eqref{coercivity} and the definition
\eqref{lambdaone} of $\lambda_1(h)$, a.e.,
\[
c\int_\Omega \!\left( \sum_{|\alpha|=m} \!\!\!\!\!\!\! \quad
a_\alpha(x,u,\dots , D^m_xu) D_x^\alpha u + a(x) |u|^{q+1} \right)
dx \geq \lambda_1(||u(.,t)||^2_{L^2(\Omega)}),
\]
where $c=\max(C,1)$, $C$ from \eqref{coercivity}. As a consequence, a.e.,
\[
\frac{1}{2} \frac{d}{dt} (||u(.,t)||^2_{L^2(\Omega)}) +
c\lambda_1(||u(.,t)||^2_{L^2(\Omega)}) \leq 0,\quad c>0.
\]
We have an ordinary differential inequality for the function
$\displaystyle y(t) = ||u(.,t)||^2_{L^2(\Omega)}$. Therefore the end of the proof is straightforward by solving of obtained differential inequality. $\fin$\\
\vspace{0 cm}\\
Now, from Proposition \ref{propositionkeystone}, we need an
estimate for $\lambda_1(h)$ from below. For this purpose, rough
estimates of $v_h$ in $L^\infty$-norm and $\lambda_1(h)$ by above
are indispensable. But \eqref{integralassumption} (or
\eqref{integralassumptionorlitz}) does not give directly an
a-priori estimate of $\lambda_1(h)$. It is why we use a trick.
Suppose that $O$ belongs to $\Omega$. We define
\begin{equation} \label{atilde}
\widetilde{a}(x):=a(x) \exp\left(-\frac{1}{|x|^\alpha}\right),
\alpha>0.
\end{equation}
In a same way, $\displaystyle \widetilde{\lambda}_1(h) = \inf \left\{ \int_\Omega (|D^{m}_{x}v|^2 + \widetilde{a}(x) |v|^{1+q}) \, dx, v \in W^{m,2}_0(\Omega), ||v||^2_{L^2(\Omega)}=h \right\}$. Since $0 \leq \widetilde{a}(x) \leq a(x)$, $\widetilde{\lambda}_1(h) \leq \lambda_1(h)$ for all $h>0$. Hence, $\displaystyle \int_0^1 \frac{1}{\lambda_1(h)} dh \leq \int_0^1 \frac{1}{\widetilde{\lambda}_1(h)} dh$. So, if $\displaystyle \int_0^1 \frac{1}{\widetilde{\lambda}_1(h)} dh < +\infty$, we get a finite extinction time.\\
For $N \neq 2m$, if $\alpha>0$ is small enough, $\displaystyle x \mapsto \exp\left(-\frac{1}{|x|^\alpha}\right)$ satisfies
\eqref{integralassumption} by Proposition \ref{propositionexpinsvarphipower}. By Theorem \ref{theoremstabilityofstheta}, both $a(x)$ and $\widetilde{a}(x)$ satisfy the same condition \eqref{integralassumption} but $\widetilde{a}(x)$ holds the a-priori estimate
\begin{equation} \label{aprioriestimatefortildea}
\widetilde{a}(x) \leq C \, \exp\left(-\frac{1}{|x|^\alpha}\right).
\end{equation}
For $N=2m$, by Proposition \ref{propositionexpinsvarphi}, if
$\alpha>0$ is small enough, $\displaystyle x \mapsto
\exp\left(-\frac{1}{|x|^\alpha}\right)$ satisfies
\eqref{integralassumptionorlitz}. In a same way, by Theorem
\ref{theoremstabilityofstheta}, both $a(x)$ and $\widetilde{a}(x)$
satisfy the same condition \eqref{integralassumptionorlitz} but
$\widetilde{a}(x)$ holds also estimate
\eqref{aprioriestimatefortildea}. With estimate
\eqref{aprioriestimatefortildea}, we get
\begin{lemma} \label{lemmaaprioriestimateforlambdaone}
There exists some $C>0$ such that for $h>0$ small enough, \eqref{aprioriestimatefortildea} implies
\begin{equation} \label{lambdaonelemmaupperboundimproved}
\widetilde{\lambda}_1(h) \leq  C \, h \, (-\ln h)^\frac{2m}{\alpha}.
\end{equation}
\end{lemma}
\proof The proof is an adaptation of \cite{Be04}. Let $v \in C^\infty_0(B)$ ($B$ is the unit-ball of $\R^N$) with $v \geq 0$ and $||v||_{L^2(\Omega)}=1$, so by homogeneity,
\[
\widetilde{\lambda}_1(h) \leq h \int_\Omega |D^m v|^2 dx +
h^\frac{1+q}{2} \int_\Omega \widetilde{a}(x) |v(x)|^{1+q} \, dx.
\]
Let $0 < r \leq r_0$. We set $\displaystyle v_r(x)=v\left(\frac{x}{r}\right)$. Then,
\[
\int_\Omega v_r^2(x) \, dx = \int_{B_r} v_r^2(x) \, dx = \int_{B_r} v^2\left(\frac{x}{r}\right) \, dx = r^N \, \int_B v^2(y) \, dy = r^N,
\]
with a translation. As a consequence, $\displaystyle \left\| \frac{v_r}{r^\frac{N}{2}} \right\|_{L^2(\Omega)} = 1$. On the other hand,\\
$\displaystyle D^\alpha_x v_r(x) = r^{-|\alpha|} D^\alpha_\xi v\left(\frac{x}{r}\right)$. As a consequence, there exists $C>0$ such that
\[
\int_\Omega |D^m_x v_r|^2\, dx = \int_{B_r} \vert D^m_xv_r\vert^2
\, dx \leq r^N \frac{C}{r^{2m}}.
\]
Then by using $\displaystyle \frac{v_r}{r^\frac{N}{2}}$ in the definition of $\widetilde{\lambda}_1(h)$,
\[
\widetilde{\lambda}_1(h) \leq C \, \frac{h}{r^{2m}} + h^\frac{1+q}{2} \, \frac{r^N}{r^{N\frac{(1+q)}{2}}} \int_B \widetilde{a}(ry) |v(y)|^{1+q} \, dy.
\]
If we estimate $\displaystyle \int_B \widetilde{a}(ry) |v(y)|^{1+q} \, dy$ by $\displaystyle C \, \exp\left(-\frac{1}{r^\alpha}\right) \, \int_B |v(y)|^{1+q} \, dy$, then
\[
\widetilde{\lambda}_1(h) \leq C' \, \left( \frac{h}{r^{2m}} + h^\frac{1+q}{2} \, r^{N\frac{(1-q)}{2}} \exp\left(-\frac{1}{r^\alpha}\right) \right).
\]
To balance both terms, we choose $\displaystyle r=\frac{1}{\left(- \ln h\right)^\frac{1}{\alpha}}$. By substituting $r$,
\[
\widetilde{\lambda}_1(h) \leq C \, \left( h \, \left(- \ln h\right)^\frac{2m}{\alpha} + h^{\frac{1+q}{2}+1} \frac{1}{\left(- \ln h\right)^\frac{N(1-q)}{2\alpha}} \right) \leq C' \, h \, \left(- \ln h\right)^\frac{2m}{\alpha},
\]
for $h$ small enough which completes the proof. \fin\\
\vspace{0 cm}\\
We introduce the functional
\begin{equation} \label{functional}
\widetilde{F}(v)=\int_\Omega (|D^m v|^2  + \widetilde{a}(x)
|v|^{1+q}) \, dx,
\end{equation}
for all $v \in W^{m,2}_0(\Omega)$. Hence, there exists for all $h>0$, a function $\widetilde{v_h} \in W^{m,2}_0(\Omega)$ such that
\begin{equation} \label{defvh}
||\widetilde{v_h}||^2_{L^2(\Omega)}=h \; \text{and} \; \widetilde{\lambda}_1(h) \leq \widetilde{F}(\widetilde{v_h}) \leq 2 \, \widetilde{\lambda}_1(h),
\end{equation}
since $\widetilde{\lambda}_1(h) > 0$. We prove Theorem \ref{theorem} by estimating $\widetilde{\lambda}_1(h)$ from below. First, we deal with $N \neq 2m$.
\begin{proposition} \label{propositionnneqtwom}
Under assumptions \eqref{operatorl}, \eqref{monotonicity} and \eqref{integralassumption}, for $N \neq 2m$, there exist $C>0$, $C'>0$ and $\eta>0$ such that for $h$ small enough,
\begin{equation} \label{lambdaonemeasure}
C \leq \frac{\widetilde{\lambda}_1(h)}{h} \left[ \text{meas} \left\{ C' \; h^\eta \geq \widetilde{a}(x) \right\} \right]^\theta, \; \theta = \min\left(\frac{2m}{N},1\right).
\end{equation}
\end{proposition}
\proof Let $v \in W^{m,2}_0(\Omega)$ with $||v||_{L^2(\Omega)}>0$. From Definition \eqref{functional} of functional $\widetilde{F}$, it follows
\begin{equation} \label{Hnotation}
\int_\Omega |D^m v|^2 dx \, dx = \int_{ \{x : |v| > 0 \} } |v|^2
H(v,x) \, dx, \; H(v,x) =
\frac{\widetilde{F}(v)}{||v||_{L^2(\Omega)}^2} -
\frac{\widetilde{a}(x)}{|v|^{1-q}},
\end{equation}
which yields
\begin{equation} \label{beforeimbedding}
C_1 \, ||D_x^m v||^2_{L^2(\Omega)} \leq \int_{ \{x : |v| > 0 \} }
|v|^2 H(v,x)^+ \, dx,\,\,\,H(v,x)^+:=\max(0,H(v,x)).
\end{equation}
Since $v \in W^{m,2}_0(\Omega)$, from the Sobolev imbedding, it follows :
\begin{equation} \label{imbedding}
||v||_{L^{p^*}(\Omega)}^2 \leq C_3 \, \left\| D^m_x v \right\|_{L^2(\Omega)}^2,
\end{equation}
where constant $C_3$ does not depend on $v$ and $p^*$ is defined by
\begin{equation} \label{pstar}
p^* = \left\{ \begin{array}{lll} \frac{2N}{N-2m} & \text{if} & N>2m\\ +\infty & \text{if} & N<2m \end{array} \right. .
\end{equation}
Combining estimate \eqref{imbedding} and equality
\eqref{Hnotation} we obtain:
\[
C_4 \, ||v||_{L^{p^*}(\Omega)}^2 \leq \int_{ \{x : |v| > 0 \} } |v|^2 H(v,x)^+ \, dx,
\]
Using H\"older's inequality for estimating term in right-hand side
of last inequality, we obtain
\[
C_4 \, ||v||_{L^{p^*}(\Omega)}^2 \leq ||v||_{L^{p^*}(\Omega)}^2 \left[ \int_{ \{x : |v| > 0 \} } \left( H(v,x)^+ \right)^\frac{p^*}{p^*-2} dx \right]^\frac{p^*-2}{p^*},
\]
where $\displaystyle \frac{p^*}{p^*-2} = \frac{p^*-2}{p^*} = 1$  if $p^* = +\infty$. This last inequality yields to
\[
0 < C_4 \leq \left[ \int_{ \{x : |v| > 0 \} } \left( H(v,x)^+ \right)^\frac{p^*}{p^*-2} dx \right]^\frac{p^*-2}{p^*}.
\]
where $H$ is from \eqref{Hnotation}. From this estimate follows easily
\begin{equation} \label{main}
0 < C_4 \frac{\widetilde{F}(v)}{||v||_{L^2(\Omega)}^2} \left[ \text{meas} \left( \{x : |v| >0 \} \bigcap \left\{ H(v,x) \geq 0 \right\} \right) \right]^\frac{p^*-2}{p^*}.
\end{equation}
As a consequence, we obtain
\[
C_4 \leq \frac{\widetilde{F}(v)}{||v||_{L^2(\Omega)}^2} \left[ \text{meas} \left( \{x : |v| >0 \} \bigcap \left\{ H(v,x) \geq 0 \right\} \right) \right]^\frac{2m}{N}, \; N-2m > 0,
\]
\[
C_4 \leq \frac{\widetilde{F}(v)}{||v||_{L^2(\Omega)}^2} \left[ \text{meas} \left( \{x : |v| >0 \} \bigcap \left\{ H(v,x) \geq 0 \right\} \right) \right], \; N-2m < 0.
\]
Therefore, for $v=\widetilde{v_h}$,
\begin{equation} \label{Gnotation}
C_4 \leq \frac{2\widetilde{\lambda}_1(h)}{h} \left[ \text{meas} \left\{ x : \Gamma(h,x) \geq 0 \right\} \right]^\frac{2m}{N}, \; N-2m > 0,
\end{equation}
where $\displaystyle \Gamma(h,x) = \frac{2\widetilde{\lambda}_1(h)}{h} |\widetilde{v_h}|^{1-q} - \widetilde{a}(x)$,
\[
C_4 \leq \frac{2\widetilde{\lambda}_1(h)}{h} \left[ \text{meas} \left\{ x : \Gamma(h,x) \geq 0 \right\} \right], \; N-2m < 0.
\]
Now, we have to estimate $|\widetilde{v}_h|^{1-q}$ from above. By definition, we know that $||\widetilde{v}_h||_{L^2}^2=h$ so for all $\varepsilon>0$,
\[
\int_\Omega \widetilde{v}_h^2 \, dx \geq \int_{\{x : \widetilde{v}_h^2 \geq \varepsilon\}} \widetilde{v}_h^2 \, dx \geq \varepsilon \; \text{meas} \{x : \widetilde{v}_h^2 \geq \varepsilon\}.
\]
By setting $\varepsilon=h^\gamma$ with $0<\gamma<1$, we get
\begin{equation} \label{Rnotation}
h^{1-\gamma} \geq \text{meas} \{x : \widetilde{v}_h^2 \geq h^\gamma\} = \text{meas} \left\{x : R(h,x) \geq 0 \right\}, R(h,x) = |\widetilde{v}_h|^{1-q} - h^\frac{\gamma(1-q)}{2}.
\end{equation}
With this inequality, it follows
\[
\text{meas} \left\{ x : \Gamma(h,x) \geq 0 \right\} = \text{meas} \left( \left\{ x : \Gamma(h,x) \geq 0 \right\} \bigcap \left\{x : R(h,x) \geq 0 \right\} \right)
\]
\[
+ \text{meas} \left( \left\{ x : \Gamma(h,x) \geq 0 \right\} \bigcap \left\{x : R(h,x) < 0 \right\} \right).
\]
But, on one hand,
\[
\text{meas} \left( \left\{ x : \Gamma(h,x) \geq 0 \right\} \bigcap \left\{x : R(h,x) \geq 0 \right\} \right) \leq \text{meas} \left\{x : R(h,x) \geq 0 \right\} \leq h^{1-\gamma},
\]
and on the other hand,
\[
\text{meas} \left( \left\{ x : \Gamma(h,x) \geq 0 \right\} \bigcap \left\{x : R(h,x) < 0 \right\} \right) \leq \text{meas} \left\{ x : \frac{2\widetilde{\lambda}_1(h)}{h} h^\frac{\gamma(1-q)}{2} \geq \widetilde{a}(x) \right\}.
\]
As a consequence, we have
\begin{equation} \label{measandvh}
\text{meas} \left\{ x : \Gamma(h,x) \geq 0 \right\} \leq h^{1-\gamma} + \text{meas} \left\{ x : \frac{2\widetilde{\lambda}_1(h)}{h} h^\frac{\gamma(1-q)}{2} \geq \widetilde{a}(x) \right\}.
\end{equation}
Hence,
\[
C_4 \leq \left( \frac{2\widetilde{\lambda}_1(h)}{h} \right)^\frac{N}{2m} \left[ h^{1-\gamma} + \text{meas} \left\{ x : \frac{2\widetilde{\lambda}_1(h)}{h} h^\frac{\gamma(1-q)}{2} \geq \widetilde{a}(x) \right\} \right], \; N-2m > 0.
\]
From \eqref{lambdaonelemmaupperboundimproved}, $\displaystyle \frac{\widetilde{\lambda}_1(h)}{h} \leq  C \, (-\ln h)^\frac{2m}{\alpha}$ which yields for $h$ small enough,
\[
\left( \frac{2\widetilde{\lambda}_1(h)}{h} \right)^\frac{N}{2m} h^{1-\gamma} \leq 3C \, h^{1-\gamma} \, (-\ln h)^\frac{N}{\alpha} \to 0,
\]
when $h \to 0$. So there exists $C_5>0$ such that for $h$ small enough,
\[
C_5 \leq \left( \frac{2\widetilde{\lambda}_1(h)}{h} \right)^\frac{N}{2m} \text{meas} \left\{ x : \frac{2\widetilde{\lambda}_1(h)}{h} h^\frac{\gamma(1-q)}{2} \geq \widetilde{a}(x) \right\}, \; N-2m > 0.
\]
Since $\gamma>0$, there exists $C'>0$ and $\eta>0$ such that for $h$ small enough,\\
$\displaystyle \frac{2\widetilde{\lambda}_1(h)}{h} h^\frac{\gamma(1-q)}{2} \leq C' \, h^\eta$. Consequently,
\begin{equation} \label{main1}
C_5 \leq \left( \frac{2\widetilde{\lambda}_1(h)}{h} \right)^\frac{N}{2m} \text{meas} \left\{ x : C' \, h^\eta \geq \widetilde{a}(x) \right\}.
\end{equation}
If $N-2m < 0$, we have in a very similar way,
\begin{equation} \label{main2}
C_5 \leq \frac{2\widetilde{\lambda}_1(h)}{h} \text{meas} \left\{ x : C' \, h^\eta \geq \widetilde{a}(x) \right\},
\end{equation}
which leads to the conclusion. \fin\\
\vspace{0 cm}\\
\text{\it Proof of Theorem \ref{theorem} for $N \neq 2m$. }
Clearly, from \eqref{lambdaonemeasure},
\[
\int_0^1 \frac{dh}{\widetilde{\lambda}_1(h)} \leq \int_0^1 \frac{\text{meas} \left\{ x \, : \, C' \, h^\eta \geq \widetilde{a}(x) \right\}^\theta}{h} \, dh,
\]
where $\theta$ is from \eqref{integralassumption}. If we set $s=C' \, h^\eta$, $\displaystyle \frac{ds}{s} = \eta \frac{dh}{h}$ and so
\[
\int_0^1 \frac{dh}{\widetilde{\lambda}_1(h)} \leq \frac{1}{\eta} \, \int_0^{C'} \frac{\text{meas} \left\{ x \, : \, s \geq \widetilde{a}(x) \right\}^\theta}{s} \, ds.
\]
Since
\[
\int_0^1 \frac{dh}{\lambda_1(h)} \leq \int_0^1 \frac{dh}{\widetilde{\lambda}_1(h)},
\]
we get the conclusion thank to Proposition \ref{propositionkeystone}. \fin
\begin{proposition}
Under assumptions \eqref{operatorl}, \eqref{monotonicity} and \eqref{integralassumptionorlitz}, for $N=2m$, there exists $C>0$ such that for $h$ small enough,
\begin{equation} \label{lambdaonemeasureorlitz}
C \leq \frac{\widetilde{\lambda}_1(h)}{h} \left( \widehat{B}^{-1}\left(\text{meas} \left\{x :
|\widetilde{v}_h|^{1-q} \frac{2\widetilde{\lambda}_1(h)}{h} \geq \widetilde{a}(x) \right\} \right)^{-1}\right)^{-1}.
\end{equation}
where $\widehat{B}(s) = (s+1) \ln(s+1) - s $ is the complementary function of $\displaystyle B(t)=e^t-1-t$ in the sense of Orlicz space (see \cite{KraRut}).
\end{proposition}
\proof Let $v \in W^{m,2}_0(\Omega)$ with $||v||_{L^2(\Omega)}>0$.
We return to the functional $\widetilde{F}$ from
\eqref{functional} again. Let $\tilde v_h\in W^{m,2}_{0}(\Omega)$
is from \eqref{defvh}. Due to optimal imbedding (see [27]) the
following estimate holds:
\begin{equation} \label{orlicz}
||\widetilde{v}_h||_{L_A(\Omega)} \leq C_3 \, ||D_x^m\widetilde{v}_h||^2_{L^2(\Omega)},
\end{equation}
where $L_A(\Omega)$ is the Orlicz space related to $\displaystyle
A(t)=\exp\left(t^\frac{p}{p-1}\right)$ (see \cite{KraRut}) and
$C_3$ is a positive constant which does not depend on
$\widetilde{v}_h$. Thus, we deduce from \eqref{orlicz} and
\eqref{beforeimbedding} for $v=\widetilde{v}_h$ :
\[
C_4 \, ||\widetilde{v}_h||_{L_A(\Omega)}^2 \leq \int_{ \{x : |\widetilde{v}_h| > 0 \} } |\widetilde{v}_h|^2 \left( \frac{\Gamma(h,x)}{|\widetilde{v}_h|^{1-q}} \right) dx,
\]
where $\Gamma(h,x)$ is from \eqref{Gnotation}. So,
\[
C_4 \, ||\widetilde{v}_h||_{L_A(\Omega)}^2 \leq \int_{ \{x : |\widetilde{v}_h| >0 \} } |\widetilde{v}_h|^2 \left( \frac{\Gamma(h,x)}{|\widetilde{v}_h|^{1-q}} \right)^+ dx.
\]
By setting $B(t)=e^t-1-t$ and using the generalized version of H\"older's inequality \eqref{gvhi},
\[
C_4 \, ||\widetilde{v}_h||_{L_A(\Omega)}^2 \leq ||\widetilde{v}_h^2||_{L_B(\{x : |\widetilde{v}_h| > 0 \})} \, \left\| \left( \frac{\Gamma(h,x)}{|\widetilde{v}_h|^{1-q}} \right)^+ \right\|_{L_{\widehat{B}(\{x : |\widetilde{v}_h| > 0 \})}}.
\]
By Proposition \ref{propositionsquareroot}, $\displaystyle ||\widetilde{v}_h||_{L_A(\Omega)}^2 = ||\widetilde{v}_h^2||_{L_M(\Omega)}$. But $B(t) = e^t-1-t \leq e^t-1 = A(\sqrt{t})=M(t)$ implies by Proposition \ref{blessthana}, $\displaystyle ||\widetilde{v}_h^2||_{L_B(\Omega)} \leq ||\widetilde{v}_h^2||_{L_M(\Omega)}$ and so,
\[
C_4 \, ||\widetilde{v}_h^2||_{L_B(\Omega)} \leq ||\widetilde{v}_h^2||_{L_B(\{x : |\widetilde{v}_h| > 0 \})} \, \left\| \left( \frac{\Gamma(h,x)}{|\widetilde{v}_h|^{1-q}} \right)^+ \right\|_{L_{\widehat{B}(\{x : |\widetilde{v}_h| > 0 \})}}.
\]
Furthermore, $\displaystyle ||\widetilde{v}_h^2||_{L_B(\{x : |\widetilde{v}_h| > 0 \})} \leq ||\widetilde{v}_h^2||_{L_B(\Omega)}$ and as a consequence,
\[
C_4 \leq \left\| \left( \frac{\Gamma(h,x)}{|\widetilde{v}_h|^{1-q}} \right)^+ \right\|_{L_{\widehat{B}(\{x : |\widetilde{v}_h| > 0 \})}}.
\]
We have
\[
\left\| \left( \frac{\Gamma(h,x)}{|\widetilde{v}_h|^{1-q}} \right)^+ \right\|_{L_{\widehat{B}(\{x : |\widetilde{v}_h| > 0 \})}} = \left\| \left( \frac{\Gamma(h,x)}{|\widetilde{v}_h|^{1-q}} \right)^+ \right\|_{L_{\widehat{B}\left(\{x : |\widetilde{v}_h| > 0 \} \bigcap \left\{x : \Gamma(h,x) \geq 0 \right\} \right)}}
\]
\[
\leq \left\| \left( \frac{\Gamma(h,x)}{|\widetilde{v}_h|^{1-q}} \right)^+ \right\|_{L_{\widehat{B}\left( \left\{x : \Gamma(h,x) \geq 0 \right\} \right)}},
\]
by Proposition \ref{esubsetf}. With Proposition \ref{lbinfty}, we get
\[
C_4 \leq \left\| \left( \frac{\Gamma(h,x)}{|\widetilde{v}_h|^{1-q}} \right)^+ \right\|_{L^\infty(E)} \left( \widehat{B}^{-1} \left( \left( \text{meas} \left\{x : \Gamma(h,x) \geq 0 \right\} \right)^{-1} \right) \right)^{-1}.
\]
We have
\[
\left\| \left( \frac{\Gamma(h,x)}{|\widetilde{v}_h|^{1-q}} \right)^+ \right\|_{L^\infty(E)} \leq \frac{\widetilde{\lambda}_1(h)}{h},
\]
so
\[
C_4 \leq \frac{2\widetilde{\lambda}_1(h)}{h} \left( \widehat{B}^{-1} \left(\left(\text{meas} \left\{x : \Gamma(h,x) \geq 0 \right\} \right)^{-1} \right) \right)^{-1},
\]
when $h$ is small enough which leads to the conclusion. \fin\\
\vspace{0 cm}\\
\text{\it Proof of Theorem \ref{theorem} for $N = 2m$ .} If $b$ is
the derivative of $B$ then $b^{-1}(s)=\ln (s+1)$ and
\[
\widehat{B}(s)=\int_0^s b^{-1}(\sigma) \, d\sigma \sim s \ln s,
\]
when $s \to +\infty$. So there exists $C_0>0$ such that for $s$ large enough,\\
$\displaystyle \widehat{B}(s) \leq C_0 \, s \ln s = D(s)$. Hence by Proposition \ref{fg}, $\displaystyle \widehat{B}^{-1}(s) \geq D^{-1}(s)$ always for $s$ large enough. Moreover, estimate $\displaystyle \ln D(s) = \ln C_0 + \ln s + \ln(\ln s) \sim \ln s$ gives\\
$\displaystyle s \sim \frac{D(s)}{C_0 \, \ln D(s)}$, i.e., $\displaystyle D^{-1}(s) \sim \frac{s}{C_0 \, \ln s}$. So there exists some positive $K$ such that for $s$ large enough, $\displaystyle \widehat{B}^{-1}(s) \geq D^{-1}(s) \geq K \frac{s}{\ln s}$. From \eqref{lambdaonemeasureorlitz}, for $h$ small enough,
\[
C \leq \frac{\widetilde{\lambda}_1(h)}{h} \left( -\ln \text{meas} \left\{x : \Gamma(h,x) \geq 0 \right\} \right) \left( \text{meas} \left\{x : \Gamma(h,x) \geq 0 \right\} \right).
\]
If $0<\gamma<1$ then estimate \eqref{measandvh} is true, i.e.,
\[
\text{meas} \left\{x : \Gamma(h,x) \geq 0 \right\} \leq h^{1-\gamma} + \text{meas} \left\{ x : \frac{2\widetilde{\lambda}_1(h)}{h} h^\frac{\gamma(1-q)}{2} \geq \widetilde{a}(x) \right\},
\]
which implies together with estimate \eqref{lambdaonelemmaupperboundimproved} that $\displaystyle \text{meas} \left\{x : \Gamma(h,x) \geq 0 \right\} \to 0$ when $h \to 0$. If for all positive $s$, we set
\begin{equation} \label{definitionE}
E(s)=s(-\ln s),
\end{equation}
then $\displaystyle C\leq \frac{\widetilde{\lambda}_1(h)}{h} E\left( \text{meas} \left\{x : \Gamma(h,x) \geq 0 \right\} \right)$.

\noindent The function $E$ is increasing in a neighbourhood of
zero so,
\[
E^{-1}\left(\frac{Ch}{\widetilde{\lambda}_1(h)}\right) \leq \text{meas} \left\{x : \Gamma(h,x) \geq 0 \right\}.
\]
By \eqref{measandvh},
\[
E^{-1}\left(\frac{Ch}{\widetilde{\lambda}_1(h)}\right) \leq h^{1-\gamma} + \text{meas} \left\{ x : \frac{2\widetilde{\lambda}_1(h)}{h} h^\frac{\gamma(1-q)}{2} \geq \widetilde{a}(x) \right\},
\]
i.e.,
\[
1 \leq \frac{h^{1-\gamma}}{E^{-1}\left( Ch(\widetilde{\lambda}_1(h))^{-1} \right)} + \frac{\text{meas} \left\{ x : \frac{2\widetilde{\lambda}_1(h)}{h} h^\frac{\gamma(1-q)}{2} \geq \widetilde{a}(x) \right\}}{E^{-1}\left( Ch(\widetilde{\lambda}_1(h))^{-1} \right)}.
\]
But from \eqref{lambdaonelemmaupperboundimproved},
\[
\frac{h^{1-\gamma}}{E^{-1}\left( Ch(\widetilde{\lambda}_1(h))^{-1} \right)} \leq \frac{h^{1-\gamma}}{E^{-1}\left( C'' \left( -\ln h \right)^\frac{-2m}{\alpha} \right)} \to 0,
\]
since when $s \to 0$, $\displaystyle E^{-1}(s) \sim \frac{s}{-\ln s}$. Consequently, for $h$ small enough,
\[
E^{-1}\left( Ch(\widetilde{\lambda}_1(h))^{-1} \right) \leq 2 \; \text{meas} \left\{ x : \frac{2\widetilde{\lambda}_1(h)}{h} h^\frac{\gamma(1-q)}{2} \geq \widetilde{a}(x) \right\}.
\]
Always from \eqref{lambdaonelemmaupperboundimproved}, there exist $C'>0$ and $\eta>0$ such that,
\[
\text{meas} \left\{ x : \frac{2\widetilde{\lambda}_1(h)}{h} h^\frac{\gamma(1-q)}{2} \geq \widetilde{a}(x) \right\} \leq \text{meas} \left\{ x \, : \, C' \, h^\eta \geq \widetilde{a}(x) \right\},
\]
which gives
\[
C \leq \frac{\widetilde{\lambda}_1(h)}{h} E\left( 2 \; \text{meas} \left\{ x \, : \, C' \, h^\eta \geq \widetilde{a}(x) \right\} \right).
\]
We easily deduce that there exist some $K>0$ and $\delta>0$ such that
\[
K \, \int_0^\delta \frac{dh}{\widetilde{\lambda}_1(h)} \leq \int_0^\delta \left( \text{meas} \left\{ x \, : \, C' \, h^\eta \geq
\widetilde{a}(x) \right\} \right) \left( -\ln\left(\text{meas} \left\{ x \, : \, C' \, h^\eta \geq \widetilde{a}(x) \right\} \right)\right) \frac{dh}{h}.
\]
If we set $s=C' \, h^\eta$, $\displaystyle \frac{ds}{s} = \eta \frac{dh}{h}$ and so
\[
K \, \int_0^\delta \frac{dh}{\widetilde{\lambda}_1(h)} \leq \frac{1}{\eta} \, \int_0^{\delta C'} \left( \text{meas} \left\{ x \, : \, s \geq
\widetilde{a}(x) \right\} \right) \left( -\ln\left(\text{meas} \left\{ x \, : \, s \geq \widetilde{a}(x) \right\} \right)\right) \frac{ds}{s}.
\]
Since $\displaystyle \int_0^1 \frac{dh}{\lambda_1(h)} \leq \int_0^1 \frac{dh}{\widetilde{\lambda}_1(h)}$, we get the conclusion thank to Proposition \ref{propositionkeystone}. \fin\\
\vspace{0 cm}\\
We can derive some useful corollaries.
\begin{corollary}
Let $f : (0,+\infty) \to (0,+\infty)$ be a continuous
nonincreasing function such that $\displaystyle f(a(x)) \in
L^1(\Omega)$ and $\displaystyle \int_0^1 s^{-1} \, f(s)^{-\theta}
\, ds < + \infty$ where $\theta$ is defined in
\eqref{integralassumption}. Then, under assumptions
\eqref{operatorl} and \eqref{monotonicity}, for $N \neq 2m$, all
solutions of problem \eqref{parabolicequation}-\eqref{initialdata}
vanish in a finite time.
\end{corollary}
\proof If $s>0$, $\displaystyle \text{meas} \left\{ x \, : \, a(x) \leq s \right\} = \text{meas} \left\{ x \, : \, f(a(x)) \geq f(s) \right\}$ and so,
\[
\text{meas} \left\{ x \, : \, a(x) \leq s \right\} \leq f(s)^{-1} \int_\Omega f(a(x)) \, dx,
\]
and we conclude with Theorem \ref{theorem}. \fin
\begin{corollary}
Let $f : (0,+\infty) \to (0,+\infty)$ be a continuous nonincreasing function such that $\displaystyle f(a(x)) \in L^1(\Omega)$ and
\begin{eqnarray} \label{marcinkiewicz}
\int_0^1 s^{-1} \, f(s)^{-1} \, \ln f(s) \, ds < + \infty.
\end{eqnarray}
Then, under assumptions \eqref{operatorl} and \ref{monotonicity}, for $N=2m$, all solutions of problem \eqref{parabolicequation}-\eqref{initialdata} vanish in a finite time.
\end{corollary}
\proof The function $f$ has a limit when $t$ tends to zero. By \eqref{marcinkiewicz}, this limit is $+\infty$. If $s>0$ is small enough, as in the previous proof,
\[
\text{meas} \left\{ x \, : \, a(x) \leq s \right\} \leq f(s)^{-1} \int_\Omega f(a(x)) \, dx.
\]
We set $E(s)=s(-\ln s)$ for all positive $s$ and since $E$ is an increasing function in a neighbourhood of zero, there exists some $\delta>0$ such that
\[
\int_0^\delta s^{-1} E\left( M_a(s) \right) \, ds \leq \int_0^\delta s^{-1} E\left( f(s)^{-1} \int_\Omega f(a(x)) \, dx \right) \, ds,
\]
which leads to
\[
\!\!\int_0^\delta s^{-1} E\left( M_a(s) \right) \, ds \leq \!\!\left( \int_\Omega f(a(x)) \, dx \right)\!\! \int_0^\delta s^{-1} \, f(s)^{-1} \!\! \left( \ln f(s) - \ln\left( \int_\Omega f(a(x)) \, dx \right)\! \right)\! ds.
\]
But, as $f(s) \to + \infty$ when $s \to 0$, there exists some $C>0$ such that,
\[
\int_0^\delta s^{-1} E\left( M_a(s) \right) \, ds \leq C \, \int_0^\delta s^{-1} \, f(s)^{-1} \, \ln f(s) \, ds,
\]
and we conclude with Theorem \ref{theorem}. \fin\\
\vspace{0 cm}\\
There is a balance between both assumptions, i.e., $f$ has to get the right behaviour. For instance, in \cite{BHV01}, they prove that for $m=1$,
\begin{eqnarray} \label{bhvassumption}
\quad\ln \frac{1}{a} \in L^p(\Omega),
\end{eqnarray}
with $\displaystyle p > \frac{N}{2}$ implies the extinction in a finite time for the Laplacian. From the previous corollary, for more general operators,
\begin{corollary}
Under assumptions \eqref{operatorl}, \eqref{monotonicity} and \eqref{bhvassumption} for $N \neq 2m$ and $p > \theta$, all solutions of problem \eqref{parabolicequation}-\eqref{initialdata} vanish in a finite time.
\end{corollary}
We can also find a Dini-like condition in the radial case in the spirit of \cite{BS07}.
\begin{corollary}
Assume that $\displaystyle a(x) = \exp\left(-\frac{\omega(|x|)}{|x|^{N\theta}}\right)$ with $\omega$ a non decreasing and non-negative function on $(0,1]$ and $\displaystyle \omega(s) \leq \omega_0, \; \forall s \in [0,1]$. If $\omega$ satisfies
\[
\int_0^1 s^{-1}\omega(s) \, ds < + \infty,
\]
under assumptions \eqref{operatorl} and \eqref{ellipticity}, for $N \neq 2m$, one have a finite extinction time for all solutions of problem \eqref{parabolicequation}-\eqref{initialdata}.
\end{corollary}
\proof For $s>0$, $\displaystyle \text{meas} \{ x : a(x) \leq s \} = \text{meas} \left\{ x : \frac{\omega(|x|)}{|x|^{N\theta}} \geq - \ln s \right\}$. We take $x$ such that $\displaystyle \frac{\omega(|x|)}{|x|^{N\theta}} \geq - \ln s$. Since $\omega$ is bounded, $x$ satisfies\\
$\displaystyle \frac{\omega_0}{|x|^{N\theta}} \geq - \ln s$ which leads to $\displaystyle |x| \leq \left( \frac{\omega_0}{- \ln s} \right)^\frac{1}{N\theta}$. By monotonicity of $\omega$,\\
$\displaystyle \omega(|x|) \leq \omega\left(\left( \frac{\omega_0}{- \ln s} \right)^\frac{1}{N\theta}\right)$. Hence,
\[
\text{meas} \{ x : a(x) \leq s \} \leq \text{meas} \left\{ x : |x|^{N\theta} \leq \omega\left(\left( \frac{\omega_0}{- \ln s} \right)^\frac{1}{N\theta}\right) \left( - \ln s \right)^{-1} \right\}.
\]
But
\[
\!\!\text{meas} \left\{ x : |x|^{N\theta} \leq \omega\left(\left( \frac{\omega_0}{- \ln s} \right)^\frac{1}{N\theta}\right) \left( - \ln s \right)^{-1} \right\} = C_N \! \left( \!\omega\left(\left( \frac{\omega_0}{- \ln s} \right)^\frac{1}{N\theta}\!\right) \!\left( - \ln s \right)^{-1} \right)^\frac{1}{\theta}\!.
\]
So,
\[
\text{meas} \{ x : a(x) \leq s \}^\theta \leq C_N^\theta \, \omega\left(\left( \frac{\omega_0}{- \ln s} \right)^\frac{1}{\min(2m,N)} \right) \left( - \ln s \right)^{-1},
\]
which yields
\[
\int_0^\frac{1}{e} s^{-1} \, M_a(s)^\theta \, ds \leq C_N^\theta \; \int_0^\frac{1}{e} s^{-1} \, \left( - \ln s \right)^{-1} \, \omega\left(\left( \frac{\omega_0}{- \ln s} \right)^\frac{1}{N\theta}\right) \, ds.
\]
By the change of variable $\displaystyle \tau = \omega_0 \, \left( - \ln s \right)^{-1}$, that is, $\displaystyle \tau^{-1} \, d\tau = \left( - \ln s \right)^{-1} s^{-1} \, ds$,
\[
\int_0^\frac{1}{e} s^{-1} \, M_a(s)^\theta \, ds \leq C_N^\theta \; \int_0^{\omega_0} \omega\left(\tau^\frac{1}{N\theta}\right) \tau^{-1} \, d\tau.
\]
By the last change of variable $s=\tau^\frac{1}{N\theta}$, that is, $\displaystyle s^{-1} \, ds = \frac{1}{N\theta} \, \tau^{-1} \, d\tau$,
\[
\int_0^\frac{1}{e} s^{-1} \, M_a(s)^\theta \, ds \leq N\theta \, C_N^\theta \, \int_0^{\omega_0^\frac{1}{N\theta}} s^{-1} \, \omega(s) \, ds.
\]
Theorem \ref{theorem} completes the proof. \fin
\begin{corollary}
Assume that $\displaystyle a(x) = \exp\left(-\frac{\omega(|x|)}{|x|^N}\right)$ with $\omega$ a nondecreasing and nonnegative function on $(0,1]$ and $\displaystyle \omega(s) \leq \omega_0, \; \forall s \in [0,1]$. If $\omega$ satisfies
\begin{eqnarray} \label{dinicondition}
\int_0^1 s^{-1}\omega(s) \left( - \ln \left( \omega(s) \right) - \ln s \right) \, ds < + \infty,
\end{eqnarray}
under assumptions \eqref{operatorl} and \eqref{monotonicity}, for $N=2m$, one has a finite extinction time for all solutions of problem \eqref{parabolicequation}-\eqref{initialdata}.
\end{corollary}
\proof For $s>0$, $\displaystyle \text{meas} \{ x : a(x) \leq s \} = \text{meas} \left\{ x : \frac{\omega(|x|)}{|x|^N} \geq - \ln s \right\}$. We take $x$ such that $\displaystyle \frac{\omega(|x|)}{|x|^N} \geq - \ln s$. Since $\omega$ is bounded, $x$ satisfies $\displaystyle \frac{\omega_0}{|x|^N} \geq - \ln s$ which leads to\\
$\displaystyle |x| \leq \left( \frac{\omega_0}{- \ln s} \right)^\frac{1}{N}$. By monotonicity of $\omega$, $\displaystyle \omega(|x|) \leq \omega\left(\left( \frac{\omega_0}{- \ln s} \right)^\frac{1}{N}\right)$. Hence,
\[
\text{meas} \{ x : a(x) \leq s \} \leq \text{meas} \left\{ x : |x|^N \leq \omega\left(\left( \frac{\omega_0}{- \ln s} \right)^\frac{1}{N}\right) \left( - \ln s \right)^{-1} \right\}.
\]
But $\displaystyle \text{meas} \left\{ x : |x|^N \leq \omega\left(\left( \frac{w_0}{- \ln s} \right)^\frac{1}{N}\right) \left( - \ln s \right)^{-1} \right\} = C_N \, \omega\left(\left( \frac{\omega_0}{- \ln s} \right)^\frac{1}{N}\right) \left( - \ln s \right)^{-1}$.\\
So, $\displaystyle \text{meas} \{ x : a(x) \leq s \} \leq C_N \, \omega\left(\left( \frac{\omega_0}{- \ln s} \right)^\frac{1}{N}\right) \left( - \ln s \right)^{-1}$. If we set $E(s) = s(-\ln s)$ for all positive $s$, we get $\displaystyle E\left( M_a(s) \right) \leq E\left( C_N \, \omega\left(\left( \frac{\omega_0}{- \ln s} \right)^\frac{1}{N}\right) \left( - \ln s \right)^{-1} \right)$ for $s$ small enough since $E$ is an increasing function in a neighbourhood of zero. As a consequence, there exists some $\delta>0$ such that\\
$\displaystyle \int_0^\delta s^{-1} \, E\left( M_a(s) \right) \, ds \leq \int_0^\delta s^{-1} E\left( C_N \, \omega\left(\left( \frac{\omega_0}{- \ln s} \right)^\frac{1}{N}\right) \left( - \ln s \right)^{-1} \right) \, ds$,\\
which gives $\displaystyle \int_0^\delta s^{-1}E\left( M_a(s) \right) \, ds$
\[
\!\!\leq \!\int_0^\delta \frac{C_N}{s} \omega\left(\left( \frac{\omega_0}{- \ln s} \right)^\frac{1}{N}\!\!\right) \left( - \ln s \right)^{-1} \!\!\left( \!- \ln C_N - \ln \left( \omega\left(\left( \frac{\omega_0}{- \ln s} \right)^\frac{1}{N}\!\!\right) \right) + \ln \left(\!- \ln s \right) \!\right) ds.
\]
But $\omega$ satisfies \eqref{dinicondition} which means that by monotonicity, $\omega(s) \to 0$ when $s \to 0$. So for $s$ small enough, $\displaystyle - \ln C_N - \ln \left( \omega\left(\left( \frac{\omega_0}{- \ln s} \right)^\frac{1}{N}\right) \right) + \ln \left(- \ln s \right)$
\[
\leq 2 \; \left( - \ln \left( \omega\left(\left( \frac{\omega_0}{- \ln s} \right)^\frac{1}{N}\right) \right) + \ln \left(- \ln s\right) \right).
\]

Consequently, for some $0 < \delta' < \delta$, we get $\displaystyle \int_0^{\delta'} s^{-1}E\left( M_a(s) \right) \, ds$
\[
\leq 2 \, C_N \, \int_0^{\delta'} \omega\left(\left( \frac{\omega_0}{- \ln s} \right)^\frac{1}{N}\right)\left( - \ln \left( \omega\left(\left( \frac{\omega_0}{- \ln s} \right)^\frac{1}{N}\right) \right) + \ln \left(-\ln s\right) \right)(- s \, \ln s)^{-1} ds.
\]
By the change of variable $\displaystyle \tau = \omega_0 (- \ln s)^{-1}$, that is, $\displaystyle \tau^{-1} d\tau = (- s \, \ln s)^{-1} ds$,\\
$\displaystyle \int_0^{\delta'} s^{-1}E\left( M_a(s) \right) \, ds$
\[
\leq C_N \, \int_0^\frac{\omega_0}{- \ln \delta'} \omega\left(\tau^\frac{1}{N}\right)\left( - \ln \left( \omega\left(\tau^\frac{1}{N}\right) \right) - \ln \tau + \ln \omega_0 \right) \, \tau^{-1} d\tau.
\]
Hence, there exists $\delta''<\delta'$ such that the following estimate holds :
\[
\int_0^{\delta''} \!\! s^{-1}E\left( M_a(s) \right) \, ds \leq 3NC_N \! \int_0^\frac{\omega_0}{- \ln \delta''} \!\!\!\!\omega\left(\tau^\frac{1}{N}\right)\!\!\left( - \ln \left( \omega\left(\tau^\frac{1}{N}\right) \right) - \ln\left(\tau^\frac{1}{N}\right) \right)\!\tau^{-1} d\tau.
\]
By the last change of variable $s=\tau^\frac{1}{N}$, that is, $\displaystyle s^{-1} ds = \frac{1}{N} \, \tau^{-1} d\tau$,
\[
\int_0^{\delta''} s^{-1}E\left( M_a(s) \right) \, ds \leq 3N^2C_N \int_0^{\left( \frac{\omega_0}{- \ln \delta''} \right)^\frac{1}{N}} s^{-1} \omega(s) \left( - \ln \left( \omega(s) \right) - \ln s \right) \, ds.
\]
This time also, Theorem \ref{theorem} completes the proof. \fin
\section{Second order case}
Here we prove Theorem \ref{theorem1.2}. Our  proof is a detailed
analysis of sufficient condition of extinction of solutions
obtained in \cite{BHV01} (see condition
\eqref{sufficientconditionlaplacian} in Theorem 4.2 from
Appendix). They introduce the quantity
\begin{equation} \label{lambdaonelinear}
\lambda_{1,2}(h) = \inf \left\{ \int_\Omega \left(|\nabla v|^2 + \frac{1}{h^2} \, a(x) \, v^2\right) dx : v \in W^{1,2}_0(\Omega), \; \int_\Omega v^2 dx = 1 \right\}, \; h > 0.
\end{equation}
As in the previous section, for $\alpha>0$ small enough, changing function $a$ into
\[
\widetilde{a}(x)=a(x) \exp\left(-\frac{1}{|x|^\alpha}\right),
\]
does not change \eqref{integralassumptionlaplacian} but by defining in a very similar way
\[
\widetilde{\lambda}_{1,2}(h) = \inf \left\{\int_\Omega \left(|\nabla v|^2 + \frac{1}{h^2} \, \widetilde{a}(x) \, v^2\right) dx : v \in W^{1,2}_0(\Omega) \right\},
\]
we have the a-priori estimate by Corollary 2.23 in \cite{Be04},
\begin{equation} \label{estimateforwidetildelambdaonelinear}
\widetilde{\lambda}_{1,2}(h) \leq C \; (-\ln h)^\frac{2}{\alpha}.
\end{equation}
Since $\widetilde{\lambda}_{1,2}(h) \leq \lambda_{1,2}(h)$ and $\displaystyle t \mapsto \frac{\ln t}{t}$ is a decreasing function for $t$ large enough, condition \eqref{sufficientconditionlaplacian} from Theorem 4.2 (Appendix) is implied by
\begin{equation} \label{sufficientconditionlaplaciantilde}
\sum_{n=1}^\infty \frac{1}{\widetilde{\lambda}_{1,2}\left(\alpha_n^\frac{1-q}{2}\right)} \left( \ln \left( \widetilde{\lambda}_{1,2}\left(\alpha_n^\frac{1-q}{2} \right) \right) + \ln \left( \frac{\alpha_n}{\alpha_{n+1}} \right) + 1 \right) < +\infty,
\end{equation}
As in \cite{BHV01}, we transform condition \eqref{sufficientconditionlaplaciantilde} into a simpler form. The following theorem is an adaptation of Theorem 2.3 in \cite{BHV01}.
\begin{proposition} \label{propositionequivalencewithintegralcondition}
Condition \eqref{sufficientconditionlaplaciantilde} is equivalent to
\begin{equation} \label{integralconditionlaplaciantilde}
\int_0^1 \frac{1}{h \widetilde{\lambda}_{1,2}(h)} dh < \infty.
\end{equation}
\end{proposition}
\proof
By changing the sequence $\{\alpha_n\}$ into $\{\alpha_n^\frac{2}{1-q}\}$, \eqref{sufficientconditionlaplaciantilde} is equivalent to
\[
\sum_{n=1}^\infty \frac{1}{\widetilde{\lambda}_{1,2}\left(\alpha_n\right)} \left( \ln \left( \widetilde{\lambda}_{1,2}\left(\alpha_n\right) \right) + \ln \left( \frac{\alpha_n}{\alpha_{n+1}} \right) + 1 \right) < +\infty.
\]
Suppose that \eqref{sufficientconditionlaplaciantilde} holds. Then, it implies that $\alpha_n \to 0$ and that $\widetilde{\lambda}_{1,2}\left(\alpha_n\right) \to \infty$ as $n$ tends to infinity. Clearly, $h \mapsto \widetilde{\lambda}_{1,2}(h)$ is a nonincreasing function which means that $\{\widetilde{\lambda}_{1,2}\left( \alpha_n \right)\}$ is a nondecreasing sequence. We use estimate of Theorem 2.3 in \cite{BHV01}.
\[
\int_{\alpha_{n+1}}^{\alpha_n} \frac{1}{h \widetilde{\lambda}_{1,2}(h)} dh \leq \frac{1}{\widetilde{\lambda}_{1,2}(\alpha_n)} \ln \left( \frac{\alpha_n}{\alpha_{n+1}} \right), \; \forall n \geq 1,
\]
which yields
\[
\int_0^{\alpha_1} \frac{1}{h \widetilde{\lambda}_{1,2}(h)} dh \leq \sum_{n=1}^\infty \frac{1}{\widetilde{\lambda}_{1,2}(\alpha_n)} \ln \left( \frac{\alpha_n}{\alpha_{n+1}} \right) < +\infty.
\]
Conversely, suppose that \eqref{integralconditionlaplaciantilde} holds. We take the sequence $\alpha_n=n^{-n}$ as in \cite{BS07}. Indeed,\\
$\displaystyle \widetilde{\lambda}_{1,2}(\alpha_n) \leq C \, (n\ln n)^\frac{2}{\alpha}$ leads to $\displaystyle \ln \left( \widetilde{\lambda}_{1,2}(\alpha_n) \right) \leq C \, \ln n$ for $n$ large enough. Moreover,
\begin{equation} \label{lnratioalphan}
\ln \left( \frac{\alpha_n}{\alpha_{n+1}} \right) \sim \ln n \; \Longrightarrow \; \ln \left(\widetilde{\lambda}_{1,2}(\alpha_n) \right) \leq C \; \ln \left( \frac{\alpha_n}{\alpha_{n+1}} \right),
\end{equation}
always for $n$ large enough ($C$ is a generic positive constant). Clearly, by monotonicity,
\[
\int_{\alpha_{n+1}}^{\alpha_n} \frac{1}{h \widetilde{\lambda}_{1,2}(h)} dh \geq \frac{1}{\widetilde{\lambda}_{1,2}(\alpha_{n+1})} \ln \left( \frac{\alpha_n}{\alpha_{n+1}} \right), \; \forall n \geq 1.
\]
Hence, thanks to \eqref{lnratioalphan}, there exists $C>0$ such that for $n$ large enough,
\[
\int_{\alpha_{n+1}}^{\alpha_n} \frac{1}{h \widetilde{\lambda}_{1,2}(h)} dh \geq C \;  \frac{1}{\widetilde{\lambda}_{1,2}(\alpha_{n+1})} \ln \left( \frac{\alpha_{n+1}}{\alpha_{n+2}} \right).
\]
So we get $\displaystyle \sum_{n=1}^\infty \frac{1}{\widetilde{\lambda}_{1,2}(\alpha_n)} \ln \left( \frac{\alpha_n}{\alpha_{n+1}} \right) < +\infty$. This implies \eqref{sufficientconditionlaplaciantilde}. \fin\\
\vspace{0 cm}\\
\textit{Proof of Theorem \ref{theorem1.2}}.  For $h$ small enough,
we have the following estimate for $\widetilde{\lambda}_{1,2}(h)$
\cite{BHV01},
\[
0 < C \leq \text{meas} \{x \in \Omega : h^{-2}\widetilde{a}(x) \leq 3 \, \widetilde{\lambda}_{1,2}(h) \} \; (\widetilde{\lambda}_{1,2}(h))^\frac{N}{2}.
\]
For this estimate, they use the Leib-Thirring formula about the counting number with some properties of semi-classical analysis \cite{He89}. By \eqref{estimateforwidetildelambdaonelinear},
\[
0 < C \leq \text{meas} \{x \in \Omega : \widetilde{a}(x) \leq C \, h^2 (-\ln h)^\frac{2}{\alpha} \} \; (\widetilde{\lambda}_{1,2}(h))^\frac{N}{2}.
\]
So, for $h$ small enough, $C \, h^2 (-\ln h)^\frac{2}{\alpha} \leq h$ which gives
\[
\frac{1}{\widetilde{\lambda}_{1,2}(h)} \leq C \; \text{meas} \{x \in \Omega : \widetilde{a}(x) \leq h \}^\frac{2}{N}.
\]
As a consequence, for some $h_0>0$ small enough,
\[
\int_0^{h_0} \frac{1}{h \widetilde{\lambda}_{1,2}(h)} dh \leq C \; \int_0^{h_0} \frac{\text{meas} \{x \in \Omega : \widetilde{a}(x) \leq h \}^\frac{2}{N}}{h} dh.
\]
We conclude with the following arguments :
\[
\int_0^1 \frac{\text{meas} \{x \in \Omega : |a(x)| \leq t \}^\frac{2}{N}}{t} \, dt < \infty,
\]
implies by Theorem \ref{theoremstabilityofstheta} and Proposition \ref{propositionexpinsvarphipower}, $\displaystyle \int_0^1 \frac{\text{meas} \{x \in \Omega : |\widetilde{a}(x)| \leq t \}^\frac{2}{N}}{t} \, dt < \infty$
which yields $\displaystyle \int_0^1 \frac{1}{h \widetilde{\lambda}_{1,2}(h)} dh < \infty$ and then by proposition \ref{propositionequivalencewithintegralcondition},
\[
\sum_{n=1}^\infty \frac{1}{\widetilde{\lambda}_{1,2}\left(\alpha_n^\frac{1-q}{2}\right)} \left( \ln \left( \widetilde{\lambda}_{1,2}\left(\alpha_n^\frac{1-q}{2} \right) \right) + \ln \left( \frac{\alpha_n}{\alpha_{n+1}} \right) + 1 \right) < \infty.
\]
This last inequality means that
\[
\sum_{n=1}^\infty \frac{1}{\lambda_{1,2}\left(\alpha_n^\frac{1-q}{2}\right)} \left( \ln \left( \lambda_{1,2}\left(\alpha_n^\frac{1-q}{2} \right) \right) + \ln \left( \frac{\alpha_n}{\alpha_{n+1}} \right) + 1 \right) < +\infty,
\]
for some sequence $\{\alpha_n\}$. By Theorem 4.2 in Appendix, all solution vanish in a finite time. \fin
\begin{remark}
If we assume that $a(x)$ is greater than a positive constant in a neighbourhood of the boundary of $\Omega$ then the related Neumann problem can be reduced to the former Dirichlet problem. Indeed, the solution of the Neumann problem vanishes in a finite time in a neighbourhood of the boundary of $\Omega$ and up to a shift in time, the solution satisfies the Dirichlet boundary condition.
\end{remark}
\section{Appendix}
\subsection{The properties of classes $S_\varphi$}
Let $\varphi$ a function defined on $[0,\gamma]$ for some $\gamma>0$ which holds the following
properties :
\begin{enumerate}
    \item[1)] $\varphi(0)=0$,
    \item[2)] $\varphi$ is a nondecreasing function on $[0,\gamma]$,
    \item[3)] $\varphi(t)>0$, $\forall t \in (0,\gamma]$,
    \item[4)] there exist $C>0$ and $\gamma' \in (0,\gamma]$ such that for all $\alpha,\beta$ in $[0,\gamma']$,
\[
\varphi(\alpha+\beta) \leq C \, \left(\varphi(\alpha)+\varphi(\beta)\right).
\]
\end{enumerate}
We set
\[
S_\varphi = \left\{a \in L^\infty(\Omega) \; | \; \exists c > 0 \; : \; \int_0^c \frac{\varphi\left(\text{meas} \{x \in \Omega : |a(x)| \leq t \}\right)}{t} \, dt < + \infty \right\}.
\]
We start with some basic properties.
\begin{proposition} \mbox{} \vspace{-0.6 cm}\\ \mbox{}
\begin{enumerate}
    \item $a \in S_\varphi \; \Longleftrightarrow \; |a| \in S_\varphi$,
    \item $1 \in S_\varphi$ ($1$ stands for the constant function equal to $1$ on whole $\Omega$),
    \item if $\psi$ satisfies $(1)$, $(2)$, $(3)$ and $\varphi \leq \psi$ then $S_\varphi \supset S_\psi$,
    \item $a \in S_\varphi \; \Longleftrightarrow \; \forall \lambda \in \R^*, \; \lambda a \in S_\varphi$.
    \item $a \in S_\varphi \; \Longleftrightarrow \; \forall \kappa>0, \; |a|^\kappa \in S_\varphi$.
\end{enumerate}
\end{proposition}
\proof Let $a \in S_\varphi$ and $\lambda \in \R^*$. By the change of variable $t = |\lambda| \tau$,
\[
\int_0^c \frac{\varphi\left(\text{meas} \{x \in \Omega : |\lambda a(x)| \leq t \}\right)}{t} \, dt = \int_0^\frac{c}{|\lambda|} \frac{\varphi\left(\text{meas} \{x \in \Omega : |a(x)| \leq \tau \}\right)}{\tau} \, d\tau,
\]
which concludes the fourth assertion.\\
Let $a \in S_\varphi$ and $\kappa>0$. In a same way, by the change of variable $t = \tau^\kappa$,
\[
\int_0^c \frac{\varphi\left(\text{meas} \{x \in \Omega : |a(x)|^\kappa \leq t \}\right)}{t} \, dt = \kappa \int_0^{c^\frac{1}{\kappa}} \frac{\varphi\left(\text{meas} \{x \in \Omega : |a(x)| \leq \tau \}\right)}{\tau} \, d\tau.
\]
The proof is complete. \fin\\
\vspace{0 cm}\\
Clearly, power functions satisfy $(1)$, $(2)$, $(3)$ and $(4)$.
\begin{proposition} \label{propositionexpinsvarphipower}
For $\alpha>0$ small enough, the function $\displaystyle w(x)=\exp \left(-\frac{1}{|x|^\alpha}\right)$ belongs to $S_\varphi$ where $\varphi(x)=x^\beta$ with $\beta>0$.
\end{proposition}
\proof By direct calculations,
\[
\varphi(\text{meas} \{x \in \Omega : |w(x)| \leq t \})=\varphi(\text{meas} \{x \in \Omega : |x|^\alpha \leq (-\ln t) \})=C_N \, \frac{1}{(-\ln t)^\frac{N\beta}{\alpha}}. \fin
\]
The main property of the set $S_\varphi$ is its stability with respect to the product.
\begin{theorem} \label{theoremstabilityofstheta}
If $a$ and $b$ belong to $S_\varphi$ then $ab \in S_\varphi$.
\end{theorem}
\proof The assumption $a,b \in S_\varphi$ implies that $a(x)>0$ and $b(x)>0$ a.e. on $\Omega$ so
\[
\lim_{t \to 0} \text{meas} \{x \in \Omega : |b(x)|^2 \leq t \} + \text{meas} \left\{x \in \Omega : |a(x)|^2 \leq t \right\} = 0.
\]
Let $t>0$ small enough, i.e,
\[
\text{meas} \{x \in \Omega : |b(x)|^2 \leq t \} + \text{meas} \left\{x \in \Omega : |a(x)|^2 \leq t \right\} \leq \gamma'.
\]
Let us consider $\{x \in \Omega : |a(x)b(x)| \leq t \}$. Pick up $\eta>0$. Then we make a partition in the following way,
\[
\{x \in \Omega : |a(x)b(x)| \leq t \} = (\{x \in \Omega : |a(x)b(x)| \leq t \} \bigcap \{x \in \Omega : |a(x)| \geq \eta \})
\]
\[
\bigcup (\{x \in \Omega : |a(x)b(x)| \leq t\} \bigcap \{x \in \Omega : |a(x)| < \eta \}).
\]
For the first subset, if $x$ in $\Omega$ satisfies both conditions $|a(x)b(x)| \leq t$ and $|a(x)| \geq \eta$ then $|b(x)| \leq \frac{t}{\eta}$ which means that
\[
\{x \in \Omega : |a(x)b(x)| \leq t \} \bigcap \{x \in \Omega : |a(x)| \geq \eta \} \subset \left\{x \in \Omega : |b(x)| \leq \frac{t}{\eta} \right\}.
\]
Clearly,
\[
\{x \in \Omega : |a(x)b(x)| \leq t \} \bigcap \{x \in \Omega : |a(x)| < \eta \} \subset \left\{x \in \Omega : |a(x)| \leq \eta \right\}.
\]
As a consequence, for $\eta = \sqrt{t}$,
\[
\text{meas} \{x \in \Omega : |a(x)b(x)| \leq t \}
\]
\[
\leq \text{meas} \{x \in \Omega : |b(x)| \leq \sqrt{t} \} + \text{meas} \left\{x \in \Omega : |a(x)| \leq \sqrt{t} \right\}.
\]
So, since $\varphi$ is a nondecreasing function on $[0,\gamma]$,
\[
\varphi\left(\text{meas} \{x \in \Omega : |a(x)b(x)| \leq t \}\right)
\]
\[
\leq \varphi\left(\text{meas} \{x \in \Omega : |b(x)|^2 \leq t \} + \text{meas} \left\{x \in \Omega : |a(x)|^2 \leq t \right\}\right).
\]
But by $4)$,
\[
\varphi\left(\text{meas} \{x \in \Omega : |a(x)b(x)| \leq t \}\right) \leq
\]
\[
C \; \left[ \varphi\left(\text{meas} \{x \in \Omega : |b(x)|^2 \leq t \}\right) + \varphi\left(\text{meas} \left\{x \in \Omega : |a(x)|^2 \leq t \right\}\right) \right].
\]
By the previous proposition, $a^2$ and $b^2$ belong to $S_\varphi$ hence for some $c>0$ small enough,
\[
\int_0^c \frac{\varphi\left(\text{meas} \{x \in \Omega : |a(x)b(x)| \leq t \}\right)}{t} \, dt \leq
\]
\[
C \; \left( \int_0^c \frac{\varphi\left(\text{meas} \{x \in \Omega : |a(x)|^2 \leq t \}\right)}{t} \, dt + \int_0^c \frac{\varphi\left(\text{meas} \{x \in \Omega : |b(x)|^2 \leq t \}\right)}{t} \, dt \right).
\]
As a conclusion, $ab$ is in $S_\varphi$. \fin\\
\vspace{0 cm}\\
The next step is to find a new class of functions satisfying properties $1)$, $2)$, $3)$ and $4)$.
\begin{proposition}
Let $\varphi$ be a function defined on $[0,\gamma]$ for some $\gamma>0$ which satisfies $1)$, $2)$, $3)$ and
\begin{enumerate}
    \item[4')] $\varphi$ is a convex function on $[0,\gamma]$ with $\displaystyle \limsup_{t \to 0^+} \frac{\varphi(2t)}{\varphi(t)} < +\infty$. Then $\varphi$ satisfies $4)$.
\end{enumerate}
\end{proposition}
\proof Since $\varphi$ is convex on $[0,\gamma]$, for all $\alpha,\beta$ in $\displaystyle \left[0,\frac{\gamma}{2}\right]$, $\displaystyle \varphi(\alpha+\beta) \leq \frac{\varphi(2\alpha)+\varphi(2\beta)}{2}$. It remains to prove that for all $t>0$ some enough, $\varphi(2t) \leq C \; \varphi(t)$ for some $C>0$.\\
Always by convexity of $\varphi$, the function $\displaystyle t \mapsto \frac{\varphi(2t)}{\varphi(t)}$ is continuous on $\displaystyle \left(0,\frac{\gamma}{2}\right]$ and bounded in a neighbourhood of zero (this function is nonnegative).\\
As a consequence, it is bounded on $\displaystyle \left(0,\frac{\gamma}{2}\right]$. \fin
\begin{proposition} \label{propositionexpinsvarphi}
The function $\displaystyle \varphi(t) = t(-\ln t)$ satisfies $1)$, $2)$, $3)$ and $4')$ for $\displaystyle \gamma = e^{-1}$. Moreover, the function $\displaystyle w(x) = \exp \left(-\frac{1}{|x|^\frac{N}{2}}\right)$ belongs to $S_\varphi$.
\end{proposition}
\proof $1)$, $2)$, $3)$ are clear. For all $\displaystyle t \in (0,e^{-1}]$, $\displaystyle \varphi'(t) = \frac{1}{(-\ln t)}+\frac{1}{(-\ln t)^2}$ so $\varphi'$ is an increasing function, hence, $\varphi$ is strictly convex. Clearly, $\displaystyle \lim_{t \to 0^+} \frac{\varphi(2t)}{\varphi(t)} = 2$.\\
With the estimate $\varphi(t) \leq t$ for $t \in [0,e^{-1}]$, we have
\[
\int_0^{e^{-1}} \frac{\varphi\left(\text{meas} \{x \in \Omega : |w(x)| \leq t \}\right)}{t} \, dt \leq \int_0^{e^{-1}} \frac{\text{meas} \{x \in \Omega : |w(x)| \leq t \}}{t} \, dt
\]
\[
= \int_0^{e^{-1}} \frac{\text{meas} \{x \in \Omega : |x|^\frac{N}{2} \leq (-\ln t)^{-1} \}}{t} \, dt = C_N \int_0^{e^{-1}} \frac{1}{t(- \ln t)^2} \, dt < + \infty. \fin
\]
\subsection{Orlicz spaces}
Let $A$ be an $N$-function \cite{KraRut}. When the derivative of $A$ called $\overline{a}$ is increasing, the $N$-functions $A$ and $\widehat{A}$ given by
\[
A(t)=\int_0^t \overline{a}(\tau) \, d\tau, \quad \widehat{A}(t)=\int_0^t \overline{a}^{-1}(\tau) \, d\tau,
\]
are said to be complementary. The Orlicz space connected to $A$ is denoted by $L_A(\Omega)$. If $E$ is a measurable set of positive measure, the Luxemburg norm is
\[
||u||_{L_A(E)} = \inf \left\{k>0 \; : \; \int_E A\left(\frac{|u(x)|}{k}\right) dx \leq 1 \right\},
\]
if the previous set is not empty and also we have a generalized version of H\"older's inequality,
\begin{eqnarray} \label{gvhi}
\left| \int_E u(t) \, v(t) \, dt \right| \leq 2 \, ||u||_{L_A(E)} \, ||v||_{L_{\widehat{A}}(E)}.
\end{eqnarray}
\begin{theorem}
 \cite{Tru} Let $\Omega$ be a bounded domain of $\R^N$ having the cone property and $mp=N$ where $p>1$. Set $\displaystyle A(t)=\exp\left(t^\frac{p}{p-1}\right) - 1$. Then there exists the imbedding $\displaystyle W^{m,p}(\Omega) \hookrightarrow L_A(\Omega)$.
\end{theorem}
Even if $A$ is an $N$-function, $M(t) = A(\sqrt{t})$ is not necessary an $N$-function but the quantity
\[
||u||_{L_M(E)} = \inf \left\{k>0 \; : \; \int_E A\left(\sqrt{\frac{|u(x)|}{k}}\right) dx \leq 1 \right\},
\]
is well defined for a measurable set $E$ of positive measure if $\displaystyle \int_E A\left( \sqrt{\frac{|u(x)|}{k}} \right) dx \leq 1$ for some positive $k$. With this extended notation, we have the following standart propositions :
\begin{proposition} \label{propositionsquareroot}
$\displaystyle ||v||^2_{L_A(E)} = ||v^2||_{L_M(E)}$ when the quantity in the left-hand side is well defined.
\end{proposition}
\proof From the definition, $\displaystyle ||v||^2_{L_A(E)} = \inf \left\{k>0 \; : \; \int_E A\left(\frac{|v(x)|}{k}\right) dx \leq 1 \right\}^2$. So,\\
$\displaystyle ||v||^2_{L_A(E)} = \inf \left\{k>0 \; : \; \int_E A\left(\sqrt{\frac{|v(x)|^2}{k^2}}\right) dx \leq 1 \right\}^2$ gives\\
$\displaystyle ||v||^2_{L_A(E)} = \inf \left\{k^2>0 \; : \; \int_E A\left(\sqrt{\frac{|v(x)|^2}{k^2}}\right) dx \leq 1 \right\}$ and\\ $\displaystyle ||v||^2_{L_A(E)} = ||v^2||_{L_M(E)}$. \fin
\begin{proposition} \label{blessthana}
If $B \leq A$ then $\displaystyle ||v||_{L_B(E)} \leq ||v||_{L_A(E)}$ when the quantity in the right-hand side is well defined.
\end{proposition}
\begin{proposition} \label{esubsetf}
If $E \subset F$ are two measurable sets of positive measure,\\
$\displaystyle ||v||_{L_B(E)} \leq ||v||_{L_B(F)}$ when the quantity in the right-hand side is well defined.
\end{proposition}
\begin{proposition} \label{lbinfty}
If $B$ is an $N$-function and $E$ a measurable set of positive
measure then $\displaystyle ||v||_{L_B(E)} \leq \frac{||v||_{L^\infty(E)}}{B^{-1}\left(\frac{1}{\text{meas}(E)}\right)}, \forall v \in L^\infty(E)$
\end{proposition}
\begin{proposition} \label{fg}
Let $f$ and $g$ be two increasing functions defined on a neighbourhood of $+\infty$ with $\displaystyle \lim_{x \to +\infty} f(x) = \lim_{x \to +\infty} g(x) = +\infty$. If $\displaystyle f(x) \leq g(x)$ for $x$ large enough then $\displaystyle f^{-1}(x) \geq g^{-1}(x)$ for $x$ large enough.
\end{proposition}
\subsection{Previous result for the second order case}
\textbf{\large Theorem 4.2 (\cite{BHV01})} Under assumptions
\eqref{coercivity} and \eqref{growth} with $m=1$, if there exists
a decreasing sequence $\{\alpha_n\}$ of positive real numbers such
that
\begin{equation} \label{sufficientconditionlaplacian}
\sum_{n=1}^\infty \frac{1}{\lambda_{1,2}\left(\alpha_n^\frac{1-q}{2}\right)} \left( \ln \left( \lambda_{1,2}\left(\alpha_n^\frac{1-q}{2} \right) \right) + \ln \left( \frac{\alpha_n}{\alpha_{n+1}} \right) + 1 \right) < +\infty,
\end{equation}
then any weak solution of problem \eqref{parabolicequation}-\eqref{initialdata} vanishes in a finite time.\\
\vspace{0 cm}\\
{\bf Acknowledgment.} The authors are very grateful to Laurent V\'eron for useful discussions and valuable comments. Both authors have been supported by an INTAS grant through the Project INTAS 05-1000008-7921.\\
\mbox{}

\end{document}